\newcommand{\di}{{\: \rm d}}   

\newlength{\razmak}
\newlength{\nrazmak}
\newcounter{secti}
\newcounter{item}[secti]

\newcommand{\sect}[1]{\medskip \begin{center} \refstepcounter{secti}
{\sc \thesecti.\  #1 }\end{center} \medskip}

\newcommand{\R}{\mathbb{R}}

\newcommand{\be}{\begin{equation} }

\newcommand{\ee}{\end{equation}}
\newcommand{\beq}{\begin{eqnarray} }
\newcommand{\eeq}{\end{eqnarray}}
\newcommand{\beqn}{\begin{eqnarray*} }
\newcommand{\eeqn}{\end{eqnarray*}}
\newcommand{\tends}{\rightarrow}
\newcommand{\kraj}{$\quad\Box$}

\newcommand{\Cov}{{\rm Cov\, }} 

\newcommand{\E}{{\rm E\, }} 

\documentclass{amsart}

\begin{document}

\begin{center}

{\LARGE\bf
A mean value theorem for systems of integrals

}

\vspace{3em}

{\sc Slobodanka Jankovi\'c, Milan Merkle}

\bigskip

\parbox{25cc}{
{\bf Abstract. }{\small More than a century ago, G. Kowalewski stated that
for each $n$ continuous functions on a compact interval $[a,b]$,
there exists an $n$-point quadrature rule (with respect to Lebesgue measure on $[a,b]$),
which is exact for given functions. Here we generalize this result to continuous
functions with an arbitrary positive and finite measure on an arbitrary interval. The proof relies
on a version of Carath\'eodory's  convex hull theorem for a continuous
curve, that we also prove in the paper. As applications, we give a representation
of the covariance for two continuous functions of a random variable, and a most general
version of Gr\"uss' inequality.

}
}

\end{center}

\renewcommand{\thefootnote}{}

\footnotetext{\scriptsize 2000 Mathematics Subject Classification. 65D32, 26A51, 26D15, 60E15

Keywords and Phrases. Quadrature rules, Carath\'eodory's  convex
hull theorem,  Covariance, Gr\"uss' inequality.

Supported by Ministry of Science and
Environmental Protection of Serbia, project number 144021.

}


\sect{Introduction and Main Results}

More than 110 years ago, G. Kowalewski  published the following result,
in a  paper entitled (in translation from
German language) "A mean value theorem for a system of $n$ integrals".

\medskip

{\bf Theorem 1.\cite{kowal1} }{\label{T1}\em Let $x_1,\ldots , x_n$ be continuous
functions in a variable $t\in [a,b]$. There
exist real numbers $t_1,\ldots, t_n$ in $[a,b]$ and  non-negative numbers $\lambda_1,\ldots,\lambda_n$, with
$\sum_{i=1}^n \lambda_i=b-a$, such that
\[ \int_a^b x_k(t)\di t = \lambda_1x_k(t_1)+\cdots +\lambda_n x_k(t_n),\quad \mbox{for each $k=1,2,\ldots n$ }.\]
}

\medskip

In \cite{kowal2}, Kowalewski generalized Theorem 1, with $\di t$ replaced with $F(t)\di t$, where $F$
is continuous and of the same sign in $(a,b)$, and with $\sum_{i=1}^n \lambda_i=\int_a^b F(t)\di t$.
It seems that these results  have not found their proper place in the literature; they were simply forgotten.
Except citations to Kowalewski's Theorem 1 in \cite{karamceb} and
\cite{mivace}, related to Gr\"uss' and Chebychev's inequalities,  we  were not able to trace any
other attempt to  use, or to generalize these results. In fact, Theorem 1 and its generalization,
presented in \cite{kowal1,kowal2}, were honestly  proved there only for $n=2$, using a theorem
(attributed to K. Weierstrass) from Hermite's course
in Analysis \cite{hermite}. Nevertheless, there is an appealing beauty,
and a potential for applications in those statements.

In this paper, we offer a generalization of Theorem 1,
for arbitrary interval $I$  (not necessarily finite), with respect to any positive finite measure, and
with functions $x_i$ that are continuous, but (if $I$ is open or infinite) not necessarily bounded.

Our main result
is the following theorem.

\medskip

{\bf Theorem 2.} {\em For an interval $I\subset \R$, let $\mu$ be a finite positive measure on
the Borel sigma-field of $I$.
 Let $x_k$, $k=1,\ldots , n$, $n\geq 1$, be continuous functions on $I$, integrable on $I$
with respect to the measure $\mu$. Then there exist points $t_1,\ldots, t_n$ in $I$,
and non-negative numbers $\lambda_1,
\ldots, \lambda_n$, with $\sum_{i=1}^n \lambda_i = \mu(I)$, such that

\[ \int_{I} x_k(t)\di \mu (t) = \sum_{i=1}^n \lambda_i x_k(t_i),\qquad k=1,\ldots, n.\]

}

\medskip

In Section 2, we prove Theorem 2 via the following version of Carath\'eodory's convex hull theorem, which also can
be of an independent interest. We show that each point in the convex hull of
a continuous curve in $\R^n$ is a convex combination of $n$ points of the curve, rather than of $n+1$ points, which
would follow from  the classical Carath\'eodory's theorem.

\medskip

{\bf Theorem 3.} {\em Let $C: t\mapsto x(t)$, $t\in I$, be a continuous curve in $\R^n$, where $I\subset \R$ is an interval,
and let $K$ be the convex hull of the  curve $C$.
Then each $v\in K$ can be represented as a convex combination of $n$ or fewer points of the curve $C$.}

\medskip

In Section 3, we discuss  Theorem 2 in the context of quadrature rules, and  in Section 4 we apply Theorem 2
to  derive a representation of  the covariance for functions of a random variable, and to obtain a most general
form of Gr\"uss' inequality.

\sect{Proofs of Theorems 2 and 3}

{\bf Proof of Theorem 3.}  According to Carath\'eodory's theorem, any point $v\in K$ can be represented as a
convex combination of at most $n+1$ points of the
curve $C$. Therefore, there exist real numbers $t_j\in I$ and $v_j\geq 0$, $0\leq j\leq n$, such that $t_0<t_1<\cdots <t_n$,
 $v_0+\cdots + v_n=1$, and
 \be
 \label{vct1}
 v=v_0 x(t_0)+v_1x(t_1)+\cdots + v_n x(t_n).
 \ee
In the sequel, we assume that all $n+1$ points $x(t_j)$ do not belong to one hyperplane, and that the numbers $v_j$ are all
positive; otherwise, one term from (\ref{vct1}) can be obviously eliminated. Denote by $p_j(x)$, $0\leq j \leq n$,
the coordinates of the vector $x\in \R^n$ with respect to the coordinate system with the origin at $v$, and with the vector
base consisted of vectors $x(t_j)-v$, $j=1,\ldots, n$ (that is, $x=\sum_{j=1}^n p_j (x) (x(t_j)-v)$). Since
\[ v_0\left( x(t_0)-v\right) = -v_1\left( x(t_1)-v\right) - \cdots - v_n\left( x(t_n)-v\right),\]
we have that $p_j\left(x(t_0)-v\right) = -v_j/v_0<0$, $j=1,\ldots, n$, i.e. the coordinates of the vector
$x(t_0)-v$ are negative.
The coordinates of vectors $x(t_j)-v$, $j=1,2,\ldots, n$ are non-negative: $p_j\left( x(t_j)-v\right) =1$ and
$p_k\left( x(t_j)-v\right)=0$ for $k\ne j$. Since the functions $t\mapsto p_j\left(x(t)-v\right)$ are continuous, the set of
points $t\geq t_0$ at which at least one of these functions reaches zero is closed, and since it is nonempty, it has the minimum.
Denoting that minimum by $\bar{t}$, we conclude that the numbers $p_j\left(x(\bar{t})-v\right)$, $j=1,\ldots, n$,
 are non-positive and at least one of them is zero. Let $p_k\left( x(\bar{t})-v\right) =0$ and $p_j\left( x(\bar{t})-v\right) \leq 0$
 for $j\ne k$. From
 \[ x(\bar{t})-v = \sum_{j=1}^{k-1} p_j\left( x(\bar{t})-v\right)(x(t_j)-v) +
 \sum_{j=k+1}^{n} p_j\left( x(\bar{t})-v\right)(x(t_j)-v),\]
 it follows that
 \[
 \left( 1-\sum_{j=1}^{k-1} p_j\left( x(\bar{t})-v\right) -\sum_{j=k+1}^{n} p_j\left( x(\bar{t})-v\right) \right) v=\]
 \[x(\bar{t}) -\sum_{j=1}^{k-1} p_j\left( x(\bar{t})-v\right) x(t_j) - \sum_{j=k+1}^{n} p_j\left( x(\bar{t})-v\right) x(t_j),\]
 wherefrom it follows that $v$ is a convex combination of  points $x(\bar{t})$ and $x(t_j)$, $j=1,\ldots, n$, $j\ne k$.

\medskip

{\bf Proof of  Theorem 2.} For given continuous and $\mu$-integrable functions $x_k$, $k=1,\ldots, n$, defined
on an interval $I$, let $C$ be the curve in $\R^n$   parametrized with $x_1=x_1(t),\ldots, x_n=x_n(t)$, $t\in I$.
Without loss of generality, we prove the statement of Theorem 2 under
 the following two assumptions.

\begin{itemize}
\item[A1] The measure $\mu$ is probability measure, i.e, $\mu (I)=1$;
\item[A2] The functions $x_1,\ldots, x_n$ are not $\mu$-a.e affine dependent. That is,
for any  hyperplane $\pi :  \alpha_1 x_1+\cdots +\alpha_n x_n=\beta$, with $\sum_{i=1}^n \alpha_i^2>0$,
$\mu\{ t\in I\; |\; x(t) \not\in \pi\}>0$.

\end{itemize}

If the assumption A2 does not hold for the given set of functions $x_1,\ldots, x_n$, it suffices to prove
the theorem for the maximal subset, say $x_1,\ldots, x_{n_1}$, $n_1<n$, of functions that satisfy A2,
and then it follows automatically for the remaining ones, which can be expressed $\mu-$a.e.
as  affine combinations of $x_1,\ldots, x_{n_1}$. In particular case $n=1$, the condition A2  implies that the
function $x_1$ is not a constant $\mu$-a.e; if it is, the theorem is trivially true.

Let

\be
\label{intjk}
  J:=(J_1,\ldots,J_n), \quad \mbox{where} \quad J_k = \int_{I} x_k(s) \di \mu (s),\ k=1,\ldots, n.
 \ee

 Let $K$ be the convex hull of $C$.  The point $J\in \R^n$,
 defined  in (\ref{intjk}) belongs to the closure $\bar{K}$. To prove that,
 let us first suppose that $I$ is a finite closed interval $[a,b]$.
 For a positive integer $m$, let $s^{(m)}_i$, $i=0,\ldots , m$ be points in the interval $[a,b]$,
such that $a=s_0<s^{(m)}_1<\cdots < s^{(m)}_m =b$, and let
$\mu^{(m)}_i =\mu ([s^{(m)}_i,s^{(m)}_{i+1}))$, $i=0,\ldots, m-2$ and $\mu^{(m)}_{m-1}=\mu([s^{(m)}_{m-1},s^{(m)}_m])$.
Let us define

\be
\label{sumjk}
  J^{(m)}:=(J_1^{(m)},\ldots,J_n^{(m)}), \quad \mbox{where}
  \quad J_k^{(m)}= \sum_{i=0}^{m-1} x_k (s^{(m)}_i)\mu^{(m)}_i ,\ k=1,\ldots, n.
 \ee


By continuity of functions $x_k$ on $[a,b]$, integrals $J_k$ are
limits as $m\tends +\infty$ of integral sums $J^{(m)}_k$, i.e.,
\[ J= \lim_{m\tends +\infty} J^{(m)}, \]
with points $s^{(m)}_i$ chosen, for example, equidistantly.
Since $\sum_{i=0}^{m-1} \mu^{(m)}_i=1$, each $J^{(m)}$ belongs
to $K$, hence, $J\in \bar{K}$.

For an open or infinite interval $I$, there exists a sequence of  closed intervals $I_p$, $p=1,2,\ldots$, such that
$I_1\subset I_2\subset \cdots$ and $\cup_p I_p =I$. Since functions $x_k$ are integrable on $I$, we have that
\be
\label{limint}
 \lim_{p\tends +\infty} \int_{I_p} x_k(t)\frac{\di \mu (t)}{\mu (I_p)} = J_k ,\qquad k=1,\ldots, n.
\ee
As the $n$-dimensional vector of integrals under the limit in (\ref{limint}) belongs to $\bar{K}$, so does
the vector $J=(J_1,\ldots, J_n)$. This proves that $J\in \bar{K}$ for arbitrary interval $I$.
Now we will show that, in fact,
$J\in\stackrel{\circ}{K}$. Indeed, if $J$ were in the boundary of the convex set $K$,
then it would have existed a hyperplane $\pi$ containing $J$, such that the points of $K$, and in particular, all
points of the curve $C$, lie in one side of $\pi$. More precisely, there would have existed real numbers $\alpha_1,
\ldots, \alpha_n$, at least one of them being non-zero, such that
\be
\label{pieq}
 \alpha_1 (x_1(t)-J_1) + \cdots + \alpha_n (x_n(t)-J_n) \geq 0,\qquad \mbox{for all $t\in I$}.
 \ee
 By the assumption A2,  the strict inequality in (\ref{pieq}) should hold on a
 subset $I'\subset I$ with $\mu(I')>0$; hence, the integral over $I$ of the left hand side in (\ref{pieq}) would
 have been  strictly positive, but it is zero. Therefore, $J$ does not lie on the boundary of $K$.

 Now, using Theorem 3, and the fact that $J\in \stackrel{\circ}{K}\subset K$,
 we conclude that $J$ can be expressed as a convex combination
 of not more than $n$ points of the curve $C$, which ends the proof.

\sect{Theorem 2 from a viewpoint of quadrature rules}

Theorem 2 claims that, given any set of continuous functions on $I$, and a finite measure $\mu$ on $I$,
there exists a (at most $n$-point) quadrature rule which is exact for those functions. As it can be seen
by inspection of the proofs
in Section 2, this quadrature rule is not unique; a point in the interior of a convex hull can be expressed as a
convex combination in infinitely many ways. This interpretation of Theorem 2 can be compared with a well known result
from \cite{karlstud}, regarding Gaussian quadratures with respect to  Chebyshev systems of functions. A brief
explanation of these terms is in order.

Real functions $x_1,\ldots , x_m$  defined on an interval $[a,b]$ are said (see \cite{karlstud}) to constitute
a Chebyshev system  on $[a,b]$ if all functions are continuous on $[a,b]$ and

\be
\label{checo}
\left | \begin{array}{llcr}
x_1(t_1) & x_1(t_2)\ldots & x_1(t_m) \\
x_2(t_1) & x_2(t_2)\ldots & x_2(t_m)\\
\rule{\razmak}{0in}\vdots    & \rule{\razmak}{0in}\vdots & \vdots\rule{\nrazmak}{0em} \\
x_m(t_1) & x_m(t_2) \ldots & x_m(t_m)
\end{array} \right |\ \ne 0
\ee
for any choice of points $t_1,\ldots,t_m \in [a,b]$ with $t_i\ne t_j$ whenever
$i\ne j$. A classical example of a Chebyshev system on any interval $[a,b]$ is furnished with functions
$x_i(t)=t^{i-1}, i=1,\ldots, m$. The condition (\ref{checo}) is equivalent to the requirement that
no $m$ points of the curve parametrized with $x_1=x_1(t),\ldots, x_m =x_m(t), t\in [a,b]$ belong to a
hyperplane which contains the origin. Another way to express (\ref{checo}) is to require that
any function of the form $g(t) = c_1 x_1(t)+\cdots +c_m x_m (t)$, $c_i \in \R$, $\sum_{i=1}^m c^2_i >0$,
must not have more than $n-1$ different zeros on $[a,b]$,

According to \cite{marowa}, for given positive and finite measure $\mu$ on  $[a,b]$, a quadrature rule
of the form
\be
\label{gauqua}
\int_{[a,b]} f(s)\di \mu (s) = \sum_{k=1}^n A_kf(t_k) +R_n(f),\qquad A_k\in \R,\quad t_k \in [a,b]
\ee
is called Gaussian with respect to  a collection  of functions $x_1,\ldots, x_{2n}$ if (\ref{gauqua})
is exact for all functions $x_i$ in place of $f$, i.e. $R_n(x_i)=0$ for $i=1,\ldots, 2n$. A quadrature rule of the
form (\ref{gauqua}) is determined by a choice of coefficients $A_k$ and points $t_k$, $k=1,\ldots, n$.

The next theorem, which can be derived from \cite[Chapter 2]{karlstud},
claims the existence and uniqueness of a Gaussian quadrature rule with respect to
a Chebyshev system of continuous functions $x_1,\ldots, x_{2n}$  on $[a,b]$.

\medskip

{\bf Theorem 4.\cite{karlstud}} {\em There exists a unique $n$-point Gaussian
quadrature rule (\ref{gauqua}) with respect to any Chebyshev system of continuous functions
$x_1,\ldots, x_{2n}$ on a finite interval $[a,b]$. Moreover,  all coefficients
 $A_1,\ldots, A_n$ are positive.}

 \medskip

There are variations and generalizations of Theorem 4 in various directions, see, for example,
\cite{marowa} or recent paper \cite{micve}.

Clearly, Theorem 4 yields a particular case of Theorem 2 if functions $x_1,\ldots, x_n$ can be complemented
with suitably chosen functions (for example $1,t,t^2,\ldots, t^n$) to make a Chebyshev system of
$2n$ functions on interval $[a,b]$. However, Theorem 2 is much more general, it is not limited to compact
intervals, it allows unbounded functions,
and does not require the condition (\ref{checo}), which is very restrictive and
difficult  to check. The price payed for the generality
is the fact that an $n$-point quadrature rule claimed in Theorem 2 is exact for
$n$ functions instead of $2n$, as in Theorem 4.

\sect{A representation of covariance and generalized Gr\"uss' inequality}

As an application of Theorem 2, we give a representation of the covariance
of random variables $f(X)$ and $g(X)$, where $f$ and $g$ are continuous functions on an interval
$I\subset \R$, and $X$ is a random variable concentrated on $I$. The idea goes back to
Karamata \cite{karamceb},
who used Theorem 1 to prove a statement of our next theorem in  particular case of uniform distribution of $X$ on
a compact interval.

As usual,  the expectation operator $E$ is defined as
\[ \E f(X)= \int_{\R} f(t) \di \mu_X(t),\]
where $f$ is a measurable function, and $\mu_X$ is a probability measure induced by $X$ on the
Borel sigma field of $\R$.  We say that  $B\subset \R$ is a support of $X$, or that $X$ is concentrated on $B$,
 if $\mu_X(B)=1$. If $X$ is concentrated on  $B$,  then
  the domain of integration
(and the domain of $f$) can be taken to be $B$.  The covariance
for random variables $U$ and $V$ is defined as
\[ \Cov (U,V)= \E (U-\E U)(V-\E V ) = \E (UV) - \E U \E V .\]

\medskip

{\bf Theorem 5. } {\em Let $X$ be a real valued random variable concentrated on an interval  $I\subseteq \R$.
  Suppose that $f$ and $g$ are continuous functions on $I$, such that $f(X)$ and $g(X)$ have finite second
order moments. Then there exist $t_1,t_2\in I$, such that
\be
 \Cov (f(X),g(X))= \frac{1}{4} \left( f(t_1)-f(t_2)\right)\left( g(t_1)-g(t_2)\right).
\ee

}

{\bf Proof. } We will use Theorem 2 with $n=2$, with functions $x_1(t)=(f(t)-\E f(X))(g(t)-\E g(X))$  and
$x_2(t)=f(t)$, and with the probability measure $\mu=\mu_X$ induced by
the random variable $X$. Using simplified notations $F=\E f(X)$ and $G=\E g(X)$, we find  that

\beqn
\Cov (f(X), g(X)) &=&\lambda (f(t_1)-F)(g(t_1)-G)
 +
(1-\lambda)(f(t_2)-F)(g(t_2)-G)\\
 F & =& \lambda f(t_1)+ (1-\lambda)f(t_2),
 \eeqn
for some $\lambda\geq 0$ and $t_1,t_2\in [a,b]$. Replacing $F$ in the first equality with the right hand side
of the second one, we conclude that
\be
\label{ado}
 \E f(X) g(X) - FG = \lambda (1-\lambda) (f(t_1)-f(t_2))(g(t_1)-g(t_2)),
\ee
for some $\lambda \in [0,1]$ and some $t_1,t_2\in I$. By continuity of $f$ and $g$, we can find another
two numbers, call them again $t_1,t_2$, so that (\ref{ado}) holds true with $\lambda=1/2$, that is,
with the maximal possible value $\lambda(1-\lambda)=1/4$.\kraj

\medskip

A well known Gr\"uss' inequality (see \cite{dsm} for a survey) can be stated in terms of an inequality for
the covariance for $f(X)$ and $g(X)$, where $X$ is a random variable. In its original formulation, it claims  that
\be
\label{grusin}
 \left|\Cov (f(X),g(X))\right| \leq \frac{1}{4} (M_f-m_f)(M_g-m_g),
 \ee
where $X$ has a uniform distribution on a compact interval $[a,b]$, $f$ and $g$ are continuous functions on $[a,b]$,
with $m_f \leq f(t)\leq M_f$ and $m_g \leq g(t) \leq M_g$ for $t\in [a,b]$. There has been a lot of related research,
and many different versions of (\ref{grusin}) are known. The following generalization of Gr\"uss' inequality
is an immediate corollary to Theorem 5.

\medskip

{\bf Theorem 6. } {\em Let $X$ be a real valued random variable concentrated on an interval
 $I\subseteq \R$. Suppose that $f$ and $g$ are continuous functions on $I$, such that
 $m_f \leq f(t)\leq M_f$ and $m_g \leq g(t) \leq M_g$ for $t\in I$. Then
\be
\label{grineq}
\left| \Cov (f(X),g(X))\right| \leq  \frac{1}{4} \left( M_f -m_f)(M_g-m_g)\right).
\ee
}

\medskip

In particular,  Theorem 6 yields the following discrete version of Gr\"uss' inequality. Let $p_n$, $n=1,2,\ldots$ be
non-negative weights with $\sum_{n=1}^{+\infty}p_n=1$, and let
$\{u_n\}$ and $\{v_n\}$ be bounded sequences of real numbers, with
\[ u\leq u_n \leq U,\quad v\leq v_n\leq V,\qquad \mbox{\rm for each $n=1,2,\ldots $}\]
Then
\be
\label{grdisc}
 \left| \sum_{n=1}^{+\infty} p_nu_nv_n -\sum_{n=1}^{+\infty} p_nu_n \sum_{n=1}^{+\infty}p_nv_n \right|
\leq \frac{1}{4}(U-u)(V-v).
\ee

To prove this inequality, it suffices to observe that there exist functions $f$ and $g$,
continuous on $\R_{+}$, such that $f(n)=u_n$ and $g(n)=v_n$, and
such that $u\leq f(t) \leq U$ and $v \leq g(t)\leq V$ for $t>0$. Then (\ref{grineq}), applied with the random
variable $X$ which takes values $n=1,2,\ldots$ with probabilities $p_n$ yields (\ref{grdisc}).

\medskip

{\bf Acknowledgements.} We wish to  thank Gradimir V. Milovanovi\'c
for a discussion related to quadrature rules and Theorem 4. Second
author acknowledges his partial affiliation to Ra\v cunarski
fakultet, Beograd, Serbia, and Instituto de Mathematica,
Universidade Federal do Rio de Janeiro, Rio de Janeiro, Brasil.

\medskip

{\sc Mathematical Institute SANU},  Knez Mihailova 35, 11000 Belgrade, Serbia

Email: bobaj@mi.sanu.ac.yu

\medskip

{{\sc Faculty of Electrical Engineering},  Bulevar Kralja Aleksandra 73, 11020 Belgrade, Serbia}

Email: emerkle@etf.bg.ac.yu
\end{document}